\documentclass[12pt,psamsfonts,leqno,oneside,letterpaper]{amsart}
\usepackage[dvips,text={6.5truein,9truein},left=1truein,top=1truein]{geometry}
\usepackage{amssymb,amsmath,amscd,enumerate}
\usepackage[pdftex]{graphicx}
\usepackage{url}

\usepackage[colorlinks,linkcolor=blue,citecolor=blue,pdfstartview=FitH]{hyperref}

\usepackage{color}

\theoremstyle{definition}
\newtheorem{para}{}[section]

\newtheorem{remark}[para]{Remark}
\newtheorem{remarks}[para]{Remarks}
\newtheorem{notation}[para]{Notation}
\newtheorem{convention}[para]{Convention}
\newtheorem{definition}[para]{Definition}
\newtheorem{definitions}[para]{Definitions}

\newcommand\Alternatives{\begin{enumerate}[(i)]}
\newcommand\EndAlternatives{\end{enumerate}}
\newcommand\Conditions{\begin{enumerate}[(1)]}
\newcommand\EndConditions{\end{enumerate}}

\theoremstyle{plain}
\newtheorem{theorem}[para]{Theorem}
\newtheorem{lemma}[para]{Lemma}
\newtheorem{proposition}[para]{Proposition}
\newtheorem{corollary}[para]{Corollary}
\newtheorem{conjecture}[para]{Conjecture}
\newtheorem{claim}[equation]{}
\newtheorem{geoconj}[para]{Geometric Conjecture}
\newtheorem{grouptheoconj}[para]{Group-Theoretic Conjecture}
\newtheorem{imptheorem}[para]{Implication Theorem}

\numberwithin{equation}{para}
\numberwithin{figure}{section}

\newcommand\Number{\begin{para}}
\newcommand\EndNumber{\end{para}}
\newcommand\Definition{\begin{definition}}
\newcommand\EndDefinition{\end{definition}}
\newcommand\Definitions{\begin{definitions}}
\newcommand\EndDefinitions{\end{definitions}}
\newcommand\GeoConj {\begin{geoconj}}
\newcommand\EndGeoConj{\end{geoconj}}
\newcommand\GroupTheoConj{\begin{grouptheoconj}}
\newcommand\EndGroupTheoConj{\end{grouptheoconj}}
\newcommand\ImpTheorem {\begin{imptheorem}}
\newcommand\EndImpTheorem{\end{imptheorem}}
\newcommand\Theorem{\begin{theorem}}
\newcommand\EndTheorem{\end{theorem}}
\newcommand\Conjecture{\begin{conjecture}}
\newcommand\EndConjecture{\end{conjecture}}
\newcommand\Remark{\begin{remark}}
\newcommand\EndRemark{\end{remark}}
\newcommand\Remarks{\begin{remarks}}
\newcommand\EndRemarks{\end{remarks}}
\newcommand\Convention{\begin{convention}}
\newcommand\EndConvention{\end{convention}}
\newcommand\Notation{\begin{notation}}
\newcommand\EndNotation{\end{notation}}
\newcommand\Lemma{\begin{lemma}}
\newcommand\EndLemma{\end{lemma}}
\newcommand\Proposition{\begin{proposition}}
\newcommand\EndProposition{\end{proposition}}
\newcommand\Corollary{\begin{corollary}}
\newcommand\EndCorollary{\end{corollary}}
\newcommand\Claim{\begin{claim}}
\newcommand\EndClaim{\end{claim}}
\newcommand\Proof{\begin{proof}}
\newcommand\EndProof{\end{proof}}
\newcommand\Equation{\begin{equation}}
\newcommand\EndEquation{\end{equation}}

\newcommand\Bullets{\begin{itemize}}
\newcommand\EndBullets{\end{itemize}}

\newcommand\vol{\mathop{\rm vol}}

\newcommand\rank{\mathop{{\rm rank}}}

\newcommand\la{\langle}
\newcommand\ra{\rangle}

\title{Hyperbolic $3$-manifolds with $k$-free fundamental group}
\author{Rosemary K. Guzman}
\begin{document}
\maketitle
\begin {abstract}
The results of Culler and Shalen for $2,3$ or $4$-free hyperbolic
$3$-manifolds are contingent on properties specific to and special
about rank two subgroups of a free group. Here we determine what
construction and algebraic information is required in order to make a
geometric statement about $M$, a closed, orientable hyperbolic
$3$-manifold with $k$-free fundamental group, for any value of $k$
greater than four. Main results are both to show what the formulation
of the general statement should be, for which Culler and Shalen's
result is a special case, and that it is true modulo a group-theoretic
conjecture. A major result is in the $k=5$ case of the geometric
statement. Specifically, we show that the required group-theoretic
conjecture is in fact true in this case, and so the proposed geometric
statement when $M$ is $5$-free is indeed a theorem. One can then use
the existence of a point $P$ and knowledge about $\pi_1(M,P)$ resulting from this theorem to attempt to improve the known lower bound on the volume of $M$, which is currently $3.44$ \cite [Theorem 1.5]{cs}.
\end{abstract}


\section{Introduction}\label{Introduction}
\medskip

The goal of this paper is to explore how the geometry of a closed, orientable hyperbolic $3$-manifold and its topological properties, especially its fundamental group, interact to provide new information about the manifold.

A hyperbolic $n$-manifold is a complete metric space that is locally isometric to the classical non-euclidean space ${\bf H}^n$ in which the sum of the angles of a triangle is less than $\pi$, or, equivalently, a complete Riemannian manifold of constant sectional curvature $-1$. Furthermore, one can express a hyperbolic $n$-manifold as the quotient of hyperbolic $n$-space modulo a discrete torsion-free group $\Gamma$ of orientation-preserving isometries, in turn $\Gamma$ is isomorphic to $\pi_1(M)$; it is this vantage point that we take in this paper. 

We will say a group $\Gamma$ is \emph{$k$-free}, where $k$ is a given positive integer, if every finitely generated subgroup of $\Gamma$ of rank less than or equal to $k$ is free. (Recall that the rank of a finitely generated group $G$ is the minimal cardinality of a generating set for $G$.)

A recurring theme here is the interplay between classical topological properties of a hyperbolic $3$-manifold and its geometric invariants, such as volume, 
and may even be regarded as a program for making the notion of Mostow Rigidity for hyperbolic $3$-manifolds explicit. The property of having $k$-free fundamental group bridges these ideas via the $\log(2k-1)$-Theorem (\cite[Main Theorem]{ACCS} combined with the Tameness Theorem  of \cite{Agol_Tameness} and \cite{Calegari_Gabai}), which uses geometric data about the manifold in regards to displacements of points under elements of $\pi_1(M)$ in ${\bf H}^3$ and forms the basis for the ideas of Section \ref{Rosemary's first section} of this paper.

One connection with topology is given by the first homology groups of $M$ with coefficients in ${\bf Z}_p$: Given an integer $k\ge 3$ and $M$ a closed, orientable, simple $3$-manifold with the property that $\dim H_1(M;{\bf Z}_2) \ge {\rm max}(3k-4,6)$, then either $\pi_1(M)$ is $k$-free or $M$ contains a closed, incompressible surface of genus at most $k-1$ which is not a fibroid \cite[Proposition 8.1]{CS_vol}.

Also, by a result of Jaco and Shalen in \cite{JaS}, any closed, orientable, hyperbolic $3$-manifold $M$ either satisfies the property that $\pi_1(M)$ is $2$-free or has a finite cover, $\widetilde {M}$, with the rank of $\pi_1(\widetilde {M})$ equal to $2$.
In this paper we are concerned with the following geometric statement: 

\GeoConj \label{rM is less than k-3 in words} If $M$ is a closed, orientable, hyperbolic $3$-manifold such that $\pi_1 (M)$ is $k$-free with $k\ge 5$, then when $\lambda = \log(2k-1)$, there exists a point $P$ in $M$ such that the 
set of all elements of $\pi_1(M,P)$ that are represented by loops of length less than $\lambda$ is contained in a subgroup of $\pi_1(M)$ of rank $\le k-3$.\EndGeoConj

In retrospect, results of Culler, Shalen, and Agol 
can be interpreted as special cases of this conjecture for the values $k=3$ and $k=4$; 
their work establishes those special cases of \ref{rM is less than k-3 in words} in \cite[Corollary 9.3]{ACS} and \cite[Theorem 1.4]{cs}. The present paper proves Conjecture \ref{rM is less than k-3 in words} for the value $k=5$, and also provides a method for showing what is required in general for the conjecture 
to hold for values of $k$ greater than five. Cases $k\le4$ have further geometric consequences than the aforementioned connections suggest at first glance --- for example, volume estimates for $M$ to be mentioned below.

Our main result will relate the Geometric Conjecture of \ref{rM is less than k-3 in words} to the following group-theoretic statement:

\GroupTheoConj \label{two} Given two rank $m$ subgroups of a free group whose intersection has rank greater than or equal to $m$, their join must have rank less than or equal to $m$ ($m\ge2$). \EndGroupTheoConj

This statement is the subject of Section \ref{Group-Theoretic Stuff} and was motivated by combining known results in the area as proved by Kent \cite{Kent}, Louder, and McReynolds \cite{LMc}. In the $k=4$ case of Conjecture \ref{rM is less than k-3 in words}, Culler and Shalen used Kent's result that if two rank-$2$ subgroups of a free group have rank-$2$ intersection, then they have a rank-$2$ join \cite{Kent}, but there were many details required to extend it to larger values of $k$.


The main result of this paper is the following theorem:

\ImpTheorem \label{implication} Group-Theoretic Conjecture \ref{two} with $m=k-2$ implies Geometric Conjecture \ref{rM is less than k-3 in words}. \EndImpTheorem

After an introduction to some terminology in Section \ref{Rosemary's first section}, Theorem 
\ref{implication} will be reformulated and proved as Theorem \ref {rM is less than k-3 again}. In the proof, we consider the action of $\Gamma$ on the sets of components of two disjoint subsets $X_i,X_j$ of a simplicial complex $K$, and using \cite[Lemma 5.12]{cs} and \cite[Lemma 5.13]{cs}, assuming the conclusion of the Theorem is false, we show that $\Gamma\le {\rm Isom}_+{\bf H}^3$ admits a simplicial action without inversions on a tree $T={\mathcal G}(X_i,X_j)$ with the property that the stabilizer in $\pi_1(M)$ of every vertex of $T$ is a locally free subgroup of $\pi_1(M)$, which is a topological impossibility.

Following a suggestion of Marc Culler and using an argument in Kent's paper \cite{Kent}, we shall establish the validity of the Group-Theoretic Conjecture \ref{two} for $m=3$; this is the topic of Section \ref{Matrices}. Hence, by Theorem \ref{implication}, Geometric Conjecture \ref{rM is less than k-3 in words} is established for the value of $k=5$, and we have the following theorem:

\Theorem \label{true for 5} Suppose $M$ is a closed,
orientable, hyperbolic $3$-manifold such that $\pi_1 (M)$ is
$k$-free with $k = 5$. Then when $\lambda = \log9$, there exists a point $P$ in $M$ such that the 
set of all elements of $\pi_1(M,P)$ that are represented by loops of length less than $\log9$ is contained in a subgroup of $\pi_1(M)$ of rank $\le 2$. \EndTheorem 
As a corollary, we state some geometric properties for particular values of $r_M$ in Section \ref {Matrices}.




As mentioned above, there has been much work done to motivate the general statement of Theorem \ref{implication}, and in those cases the geometric information was used to deduce lower bounds on the volume of $M$.
Specifically, as a special case of Theorem \ref{implication}, Culler and Shalen expressed the point $P$ of the conclusion as a $\log7$-semithick point (\cite[Theorem 1.4]{cs}), and using the existence of this point along with other consequences of $4$-freeness, they were able to show that $\vol M \ge 3.44$. In the $3$-free case, Anderson, Canary, Culler and Shalen, along with Agol, showed the existence of a point $P$ of $M$ of injectivity radius $(\log5)/2$ --- this is exactly the point $P$ described in Theorem \ref{implication} --- and used its existence to establish that $\vol M \ge 3.08$ in this case (\cite[Corollary 9.3]{ACS} and its predecessor \cite[Theorem 9.1]{ACCS}).

Note that a closed hyperbolic manifold with $k$-free fundamental group, for $k\ge2$, is in fact $k-1, k-2, k-3,\dots,2$- free. So, in particular, the results of Culler and Shalen \cite{cs} show that for a closed, orientable, hyperbolic $3$-manifold $M$ with $5$-free fundamental group, we have vol $M \ge 3.44$. The method of obtaining this bound is by finding lower bounds for both the nearby volume of the $\log7$-semithick point $P$, i.e. the volume of the $(\log7)/2$ neighborhood of $P$, and for the distant volume, i.e. the volume of the complement of this neighborhood. Therefore, a long range goal of the present work is to improve this bound with the added topological and geometric information that is gotten by virtue of the $5$-free assumption and the rank $\le 2$ subgroup described in Theorem \ref{true for 5}, with hopes that estimating the nearby and distant volumes of the given point $P$, under certain conditions, will lead to a refined lower bound on the global volume of $M$.  

I acknowledge and profoundly thank Peter Shalen, my former advisor, for his overall aegis, many helpful comments, and insights on this work which was based on my thesis. Also, I am grateful to Marc Culler and Dick Canary for their willingness to review this work and provide comment. Finally, results in Kent \cite{Kent}, Louder, and McReynolds \cite{LMc} made proving the case for $k=5$ possible.

\section{Lemma and Preliminaries} \label{Rosemary's first section}
\Definitions \label{locally finite} Suppose we are given a positive real number $\lambda > 0$ and that the subgroup $\Gamma \le {\rm Isom}_+({\bf H}^3)$ is discrete and cocompact (and so purely loxodromic). For $\gamma \in \Gamma$ we define the hyperbolic cylinder $Z_{\lambda}
(\gamma)$ to be the set of points $P \in {\bf H}^3$ such that
$d(P,\gamma \cdot P) < \lambda$. Recall that since $\gamma$ is
loxodromic, there is a $\gamma$-invariant line, $A(\gamma) \subset
{\bf H}^3$, called the \emph{axis} of $\gamma$, such that $\gamma$
acts on the points of $A(\gamma)$ as a translation by a distance $l >
0$, called the \emph{translation length} of $\gamma$. For any point $P
\in {\bf H}^3$, we have $d(P, \gamma \cdot P) \ge l$ with equality
only when $P \in A(\gamma)$. Then as long as $l < \lambda$, the
cylinder $Z_{\lambda}(\gamma)$ is non-empty (the radius of this
cylinder is computed by a simple application of the hyperbolic law of
cosines and is a monotonically increasing function for $\lambda$ in
the interval $(l, \infty)$
; see, for example, \cite{csmodp} for further details). 

\Remark \label {about a manifold}
Given $M$ a closed, orientable, hyperbolic $3$-manifold, we may
write $M$ as the quotient ${\bf H}^3 / \Gamma$, where $\Gamma$ is a
discrete group of orientation-preserving isometries of ${\bf
H}^3$ that is torsion-free. Every isometry of $\Gamma$ is
loxodromic since $M$ is closed (and so cannot be parabolic or
elliptic). Every non-trivial element $\gamma$ of $\Gamma$ is
contained in a unique maximal cyclic subgroup $C(\gamma)$ of
$\Gamma$ which is the centralizer of $\gamma$ in $\Gamma$, which means that non-trivial elements of
distinct maximal cyclic subgroups do not commute. \EndRemark

\Definition \label{cylinders} Supposing $M={\bf H}^3 / \Gamma$ is given as above, let $C(\Gamma)$ be the set of maximal cyclic subgroups of $\Gamma$. After fixing a positive real number $\lambda$, let $C_{\lambda}(\Gamma)$ denote
the set of maximal cyclic subgroups $C=C(\gamma)$ of $\Gamma$ having
at least one (loxodromic) generator $\gamma_0$ of $C$ with translation length less than $\lambda$. 
\EndDefinition

Given a cyclic subgroup $C$ of $\Gamma$, we define the set
$Z_{\lambda}(C)= \bigcup_{1 \not = \gamma \in C} Z_{\lambda}(\gamma)$. 
Then if $C \in C_{\lambda}(\Gamma)$, the set $Z_{\lambda}(C)$ is in
fact a cylinder of points in ${\bf H}^3$ that are displaced
by a distance less than $\lambda$ by some non-trivial element of $C$:
specifically, there is a loxodromic element $\gamma \in C-\{ 1 \}$ such that $Z_{\lambda} (C)=
Z_{\lambda} (\gamma)$ ($\gamma$ need not necessarily be $\gamma_0$,
the generator of $C$). Observe that if $C \in C(\Gamma) - C_{\lambda}(\Gamma)$, we have $Z_{\lambda}(C)=\emptyset$.

Note that the family of cylinders $(Z_{\lambda}
(\gamma))_{1\neq\gamma\in\Gamma}$
is \emph{locally finite} as $\Gamma$ is discrete; i.e. for every point $P$ in
${\bf H}^3$,
there is a neighborhood of $P$ which has non-empty intersection with
only finitely many of the subsets $Z_{\lambda}
(\gamma)$. Further, because the
family $(Z_{\lambda}
(\gamma))_{1\neq\gamma\in\Gamma}$ is locally finite, so then is the family $(Z_{\lambda}
(C))_{C\in C_\lambda(\Gamma)}$.\EndDefinitions

\Remark A locally finite family ${\mathcal Z}=(Z_{\lambda}
(C))_{C\in C_\lambda(\Gamma)}$ of cylinders has a natural association to the set of maximal cyclic subgroups of $\Gamma$, and if this family of cylinders covers ${\bf H}^3$, how it does so will be of particular importance, as we will later encode this information in the nerve (see Definition \ref{nerve}) of the cover ${\mathcal Z}$. Determining a ``connectedness'' argument for certain skeleta of the nerve in order to show homotopy-equivalence to ${\bf H}^3$ (and therefore contractibility), exhibited new challenges and many refinements in extending the $4$-free arguments to the $k$-free arguments and are detailed in Section \ref{Contractibility Arguments and Gamma-labeled complexes}. \EndRemark

The following lemma is an application of the $\log(2k-1)$ Theorem (\cite[Main Theorem]{ACCS} with \cite{Agol_Tameness} and \cite{Calegari_Gabai}).

\Lemma\label{the first one} Suppose $\Gamma \le Isom_+({\bf H}^3)$ is
discrete, loxodromic, $k$-free $(k\ge2)$ and torsion-free. If
there exists a point $P \in Z_{\log {(2k-1)}}(C_1) \cap\cdots\cap
Z_{\log {(2k-1)}}(C_n)$, then the rank of $\la C_1,\ldots,C_n \ra$ is $\le
k-1$. \EndLemma

\Proof (by induction on $n$)

\emph{Base case:} If $n=1$, then $P \in Z_{\log {(2k-1)}}(C)$.
Because ${\rm rk}\;C = 1$ and $k \ge 2$, ${\rm rk}\;C \le k-1$ is satisfied.

\emph{Induction assumption:} If $n=q$ then $X_q = \la C_1, \dots ,C_q \ra$, and so we assume that ${\rm rk}\;X_q \le k-1$.

\emph{Induction step:} Notice that $X_{q+1} = \la X_q, C_{q+1}\ra =
\la C_1,\dots,C_q, C_{q+1} \ra$. We must show that ${\rm rk}\;X_{q+1} \le k-1$. To simplify notation, let $r={\rm rk}\;X_q$. First, consider when ${\rm rk}\;\la X_q,C_{q+1}\ra =
r$. Since $r \le k-1$ by our induction assumption, we are done.

Next, consider the case when ${\rm rk}\;\la X_q, C_{q+1}\ra
> {\rm rk}\;X_q = r$.

\Remark As $X_q \le \Gamma$ which is $k$-free, ${\rm rk}\; X_q<k$, $C_{q+1}= \la t \ra$ is cyclic, and ${\rm rk}\;(X_q
\vee C_{q+1}) > {\rm rk}\;X_q = r$, we have $(X_q \vee C_{q+1})$ is the free product of $X_q$ and $C_{q+1}$ by
\cite[Lemma 4.3]{cs}. \EndRemark

By the remark and our induction assumption, ${\rm rk}\;\la X_q,
C_{q+1}\ra = r+1 \le (k-1)+1 = k$. Therefore ${\rm rk}\; \la X_q,
C_{q+1}\ra \le k$, leaving two subcases to consider. First, if
$r<k-1$, then ${\rm rk}\; \la X_q, C_{q+1}\ra < k$ and we are done.

In the second subcase, suppose $r=k-1$. The remark then gives that
${\rm rk}\;\la X_q, C_{q+1}\ra = r+1 = k$; we proceed to prove that ${\rm rk}\;\la X_q, C_{q+1}\ra \le k-1$ by way of contradiction.

Since $n=q+1$, by hypothesis $P \in Z_{\log {(2k-1)}}(C_1)
\cap\cdots\cap Z_{\log {(2k-1)}}(C_{q+1})$. Choose a generator
$\gamma_i$ for each $C_i$, where $1 \le i \le q+1$. For each $i$, there
exists a number $m_i \in {\bf N}$ with $d(P,{\gamma}_i^{m_i} \cdot P) <
{\rm log}(2k-1)$ by definition of the cylinders; denote this property (*).

Now the rank of $\la \gamma_1,\dots,\gamma_{q+1} \ra$ is $k$, and so this group is free (being a subgroup of $\Gamma$ which is $k$-free). In particular, $\{ \gamma_1,\dots,\gamma_{q+1} \}$ is a generating set of a free group of rank $k$, and so it must contain a subset $S$ of $k$ independent elements whose span has rank $k$. So let $S= \{ \gamma_{i_1},\dots,\gamma_{i_k} \} \subseteq \{ \gamma_1,\dots,\gamma_{q+1} \}$ be as described. Furthermore, the set $S' = \{ \gamma^{m_{i_1}}_{i_1},\dots,\gamma^{m_{i_k}}_{i_k} \}$ is also a set of $k$ independent elements whose span has rank $k$. Then as $S' \subseteq {\rm Isom}_+({\bf H^3})$ is a
set of $k$ freely-generating (loxodromic) generators with $\rank \la S' \ra = k$, 
the $\log{(2k-1)}$ Theorem of \cite{ACCS} applies here to give that
$\max_{1\le j \le k} d(P,\gamma_{i_j}^{m_{i_j}} \cdot
P) \ge \log {(2k-1)}$, thereby contradicting property (*) above.
Therefore, ${\rm rk} \; \la X_q, C_{q+1}\ra \le k-1$ as required,
and in particular is equal to $k-1$ in this subcase. \EndProof


We now provide the new notation necessary for setting up the arguments in the remaining sections, as well as a proposition relating the preceding lemma to our new notation. For the next definition, recall Definitions \ref{cylinders}. 

 \Definition \label{to G_P or not to G_P} Given a point $P \in {\bf
H}^3$, let $C_P (\lambda)$ denote the set of all $C$ in
$C_{\lambda}(\Gamma)$ for which $P$ is an element of
$Z_{\lambda}(C)$. We then associate to each point $P$ in ${\bf H}^3$ a group, $G_P(\lambda)$, which is defined by $G_P(\lambda) = \la C : C \in
C_P (\lambda) \ra $. If $C_P(\lambda) = \emptyset$, then set
$G_P(\lambda)=\la 1 \ra$, and define ${\rm rk}\;G_P(\lambda)=0$.
Also, if the value of $\lambda$ is understood to be fixed, we may
refer to $G_P(\lambda)$ simply as $G_P$.\EndDefinition

\Proposition \label{the truth about GP} Given $\Gamma \le {\rm Isom}_+({\bf
H}^3)$ is discrete, purely loxodromic, and $k$-free with $k\ge2$,
then for any point $P \in {\bf H}^3$, we have ${\rm rk }\;G_P (\log {(2k-1)})
\le k-1$. \EndProposition

\Proof This result is a direct consequence of Lemma \ref{the first one} along with the preceding definitions.
\EndProof 

\Definition \label{mer} Suppose $H$ is a subgroup of a group $G$. Then we define the \emph{minimum enveloping rank of $H$},
or $r_H$ to be the smallest rank among the ranks of groups for which
$H$ is a subgroup, if such a number exists. If $H$ is not contained in
a finitely generated subgroup of $G$, then we define $r_H$ to be
$\infty$.  More formally, when $H$ is contained in a finitely
gerenated subgroup $K$ of $G$, we may define $r_H$ as the smallest positive integer among the set $\{ \rank K : H \le K \le G \}$. \EndDefinition 

\Number \label{about mer} Note that if $H$ is non-trivial and non-cyclic, $r_H \ge 2$. Furthermore, if $h$ denotes the rank of $H$, since $H$ is in particular a subgroup of itself, by definition we have $r_H \le h$.\EndNumber

\Definition \label{the rank of a manifold} Suppose $M={\bf H}^3 / {\Gamma}$ is a closed,
orientable, hyperbolic $3$-manifold ($\Gamma \le {\rm Isom}_+({\bf H}^3)$
is discrete and purely loxodromic). 
Given a number $\lambda > 0$, we define the number $r_M(\lambda) \in {\bf N} \cup \{ 0 \}$ to be the infimum of the set $\{ r_{G_P(\lambda)} : P \in {\bf H}^3 \}$. If the value of $\lambda$ is understood to be fixed, we may refer to $r_M(\lambda)$ simply as $r_M$. \EndDefinition


Given $M={\bf H}^3/ \Gamma$ a closed, orientable, hyperbolic $3$-manifold, we now make a few observations regarding the number $r_M$:

\Number \label{about P} Given a point $P$ in ${\bf H}^3$,
it follows from the definitions that $ r_M\le r_{G_P(\lambda)}\le {\rm rk}\ G_P(\lambda)$.
\EndNumber

\Number \label{the truth about rM} When $\lambda =
\log {(2k-1)}$, as a direct consequence of Corollary \ref{the truth about
GP} and \ref{about P}, we have $r_M \le k-1$. 
\EndNumber


\Remark Notice in the standard terminology, saying that the manifold
$M$ contains a ``$\lambda$-thick'' point (i.e. a
point of injectivity radius at least $\lambda/2$ in $M$)  is reinterpreted here
as saying that $r_M(\lambda) = 0$. 
We observe that $r_M(\lambda) \not = 0$ if and only if the family of
cylinders ${\mathcal Z} = (Z_{\lambda}(C))_{C\in C_\lambda(\Gamma)}$ forms an open cover of ${\bf H}^3$. \EndRemark

When $r_M \ge 1$, we claim:
\Number \label{there is a point} 
${\bf H}^3 = {\bigcup}_{C_1,\dots,C_{r_M} \in
C_{\lambda}(\Gamma)} Z_{\lambda}(C_1) \cap \cdots \cap
Z_{\lambda}(C_{r_M})$. \EndNumber

\Proof Suppose $P$ is a point of ${\bf H}^3$. As \ref
{about P} says that ${\rm rk}\ G_P \ge r_M$, there exist maximal cyclic subgroups $C^P_1,\dots, C^P_{r_M}$ of $\Gamma$ such that $\la C_1^P, \dots,
C_{r_M}^P \ra \le G_P$ with $P \in Z_{\lambda}(C_1^P) \cap\cdots \cap
Z_{\lambda}(C_{r_M}^P)$ (keeping in mind that $P$ may be in additional
cylinders). The statement follows.
\EndProof 

\section{$\Gamma$-labeled complexes and Contractibility Arguments}\label{Contractibility Arguments and Gamma-labeled complexes}

\Definitions \label{nerve} An indexed covering ${\mathcal U}= (U_i)_{i \in
  I}$ of a topological space by non-empty open sets defines an
abstract simplicial complex called the \emph{nerve of $\mathcal U$}, denoted $K({\mathcal U})$, whose vertices are in bijective
correspondence with the elements of the index set $I$ and whose
simplices $\{ v_{i_0},\dots, v_{i_n} \}$ correspond to the non-empty
intersections $U_{i_0} \cap\cdots\cap U_{i_n}$ of sets of
${\mathcal U}$. We endow the space which is the geometric realization $|K|=|K({\mathcal U})|$ with the weak topology. Given a group $\Gamma$, a \emph{$\Gamma$-labeled
  complex} is a pair $(K, (C_v)_v)$ where $K$ is a simplicial complex and where $(C_v)_v$ is a family of cyclic subgroups of $\Gamma$ indexed by (and ranging over) the vertices $v$ of $K$.

Suppose additionally that we are given a positive real number $\lambda > 0$ and subgroup $\Gamma \le {\rm Isom}_+({\bf H}^3)$ which is discrete and cocompact. 
In particular, 
if ${\mathcal Z}(\lambda)=(Z_{\lambda}(C_i))_{i \in I, C_i \in
  C_{\lambda}(\Gamma)}$ is a cover of ${\bf H}^3$ by cylinders, then
the family ${\mathcal Z}(\lambda)$ gives rise to a $\Gamma$-labeled
  complex $(K,(C_{v_i})_{v_i})$ where $K$ is the nerve of ${\mathcal Z}(\lambda)$ and where $C_{v_i}$ is the (infinite) maximal cyclic subgroup of $\Gamma$ that corresponds to the element $Z_{\lambda}(C_{v_i})=Z_{\lambda}(C_i)$ of the cover ${\mathcal
Z}(\lambda)$ as indexed by the vertex $v_i$ of $K$. For purposes of
notation, we may refer to this vertex $v_i$ by $v_{C_i}$. \EndDefinitions

\Definition \label {labeling-compatible} Given a group $\Gamma$ and $(K, (C_v)_v)$ a
$\Gamma$-labeled complex, we say the labeling defines a
\emph{labeling-compatible} $\Gamma$-action on $(K, (C_v)_v)$ if for
every vertex $v$ of $K$, the action defined by $C_{\gamma\cdot v}
= \gamma C_v \gamma ^{-1}$ is simplicial.  \EndDefinition

\Remark Note that if $\Gamma \le {\rm Isom}_+({\bf H}^3)$ is discrete
and torsion-free, if the family ${\mathcal Z}(\lambda)  = (Z_{\lambda}(C))_{C \in
  C_{\lambda}(\Gamma)}$ covers ${\bf H}^3$, and if $K$ is the nerve of ${\mathcal Z}(\lambda)$, then the $\Gamma$-labeled complex $(K, (C_v)_v)$ admits a
labeling-compatible $\Gamma$-action. Let $V= \{ v_0, \dots, v_n \}$ be the set of vertices of an $n$-simplex of $K$; by definition $\cap_{0 \le i
\le n} Z_{\lambda}(C_{v_i}) \neq \emptyset$. Given $1 \not = \gamma \in
\Gamma$ and $v_i \in V$, define $C_{w_i} =C_{\gamma \cdot v_i}=\gamma C_{v_i} \gamma^{-1}$. We must show two things: first, $w_i$ is well-defined as a vertex of $K$; or, equivalently, that $C_{w_i}$ is a maximal cyclic subgroup in $C_{\lambda}(\Gamma)$, making the action $\gamma \cdot v_i := w_i$ a well-defined action of $\Gamma$ on the vertices of $K$; second, we must show that the set $W= \{ w_0,
\dots, w_n \}$ of vertices of $K$ is in fact the vertex set of a simplex of $K$, making this action simplicial, and therefore labeling-compatible. 
By showing that $\cap_{0 \le i \le n} Z_{\lambda}(C_{w_i})$ is non-empty, we achieve both of these goals. 

By our definition, $\cap_{0 \le i \le n}
Z_{\lambda}(C_{w_i}) = \cap_{0\le i \le n} Z_{\lambda}(C_{\gamma
\cdot v_i}) = \cap_{0 \le i \le
  n} Z_{\lambda}(\gamma C_{v_i} \gamma^{-1})$.
Using the definition of the cylinders (along with the fact that $\gamma^{-1} \in {\rm Isom}_+({\bf H}^3)$ for the first equality), 
we have $\cap_{0 \le i \le
  n} Z_{\lambda}(\gamma C_{v_i} \gamma^{-1}) = \cap_{0 \le i \le n} \gamma \cdot  Z_{\lambda}(C_{v_i}) = \gamma \cdot \cap_{0 \le i \le n} Z_{\lambda}(C_{v_i})$, which is non-empty as $\cap_{0 \le i \le n} Z_{\lambda}(C_{v_i})$ is non-empty, and hence $\cap_{0 \le i \le n}
Z_{\lambda}(C_{w_i}) \neq \emptyset$ as required. Further, since $C_{\gamma ^{-1}
\cdot w_i}= \gamma^{-1} C_{w_i} \gamma = \gamma^{-1} \gamma C_{v_i} \gamma
^{-1} \gamma = C_{v_i}$, the simplex $\{ w_0,
\dots, w_n \}$ is in fact an $n$-simplex of $K$.
\EndRemark

\Definition \label{thetasigma} Given a $\Gamma$-labeled complex
$(K,(C_v)_v)$ and $\sigma$ an open simplex in $K$, define the subgroup
$\Theta (\sigma)$ of $\Gamma$ to be the group $\la C_v : v \in
\sigma \ra$. \EndDefinition

\Number \label{thetasigma subgroup of G_P} Suppose $K$ is given to be the nerve of a family ${\mathcal Z}(\lambda) = (Z_{\lambda}(C_i))_{i \in I,
  C_i \in C_{\lambda}(\Gamma)}$ which is a cover of ${\bf H}^3$ by cylinders. If there exists a point $P \in {\bf H}^3$ in the intersection $Z_{\lambda}(C_0) \cap\cdots\cap Z_{\lambda}(C_n)$, it follows that $\{ v_{C_0},\dots,v_{C_n} \}$ is an $n$-simplex $\sigma$ of $K$, and by the Definitions \ref{to G_P or not to G_P} and \ref{thetasigma}, we have $\theta(\sigma) \le G_P(\lambda)$. \EndNumber 


\Definitions \label{mer complex} Suppose $(K,(C_v)_v)$ is a $\Gamma$-labeled complex. Given an open simplex $\sigma$ in $K$, the \emph{minimum enveloping rank of $\sigma$} will denote the minimum enveloping rank of the associated subgroup $\Theta(\sigma)$ in $\Gamma$. Notice that if $\tau \in K$ is a face of $\sigma \in K$, then we have $r_{\theta(\tau)} \le r_{\theta(\sigma)}$; i.e. the minimum enveloping rank of a face of $\sigma$ is less than or equal to that of $\sigma$. We may therefore define a subcomplex \emph{$K_{(n)}$} of $K$ to be the subcomplex that consists of the non-trivial open simplices $\sigma$ for which $r_{\theta(\sigma)} \le n$.\EndDefinitions

\Proposition \label{homoequiv} Suppose $K$ is a simplicial complex, and $\sigma$ a simplex of $K$. Suppose further that the link ${\rm lk}_K(\sigma)$ is contractible and $X \subset |K|$ is a saturated subset that contains all the simplices for which $\sigma$ is a face. Then $X - \sigma \hookrightarrow X$ is a homotopy equivalence. \EndProposition

\Proof Let $C=\cup_{\tau < \sigma} \tau$ be the union of simplices $\tau \in K$ for which $\sigma$ is a face. By how it is defined, $C$ is homotopy equivalent to ${\rm lk}_K(\sigma)$, and is therefore contractible by our assumption. 
Let $S =  {\rm star}_K\sigma$. Now $X= (X-\sigma) \cup S$ and $C= (X-\sigma) \cap S$. As $C$ and $S$ are both contractible, then by exactness of the Mayer-Veitoris sequence, the Van Kampen Theorem, and Whitehead's Theorem, it follows that $X - \sigma \hookrightarrow X$ is a homotopy equivalence.
\EndProof

\Lemma\label{rank argument general} Let $M= {\bf H}^3 /
\Gamma$. Suppose ${\mathcal Z}(\lambda) = (Z_{\lambda}(C_i))_{i \in I,
  C_i \in C_{\lambda}(\Gamma)}$ is a cover of ${\bf H}^3$ by cylinders and that $r_M \ge k-2$. Let $|K|$ denote the geometric realization of the nerve of
${\mathcal Z}(\lambda)$. Then $|K|-|K_{(k-3)}|$ is
homotopy-equivalent to ${\bf H}^3$ and therefore contractible.
\EndLemma

\Proof 
The family $(Z_{\lambda}(C_i))_{i \in I, C_i
\in C_{\lambda}(\Gamma)}$ covers ${\bf H}^3$ and has the property that every finite
intersection of (open) cylinders is contractible, as any such intersection is either
empty or convex. Thus Borsuk's Nerve Theorem \cite{bjo} applies, and we have $|K|$
is homotopy-equivalent to ${\bf H}^3$. It is only left to show that
$|K|- |K_{(k-3)}|$ is homotopy-equivalent to $|K|$. Suppose $\sigma$
is a non-trivial open simplex of $|K_{(k-3)}|$, which by definition is
to say that the minimum enveloping rank of $\theta(\sigma)$ is $\le k-3$. Let $v_{i_0}^{\sigma},..., v_{i_l}^{\sigma}$ be the vertices of $\sigma$, and set $I_{\sigma} = \{ i \in I : v_i \in \sigma \}$.

Let $U_i$ for $i \in I$ denote the cylinder $Z_{\lambda}(C_i)$ associated with the vertex $v_i$ as
defined by the nerve of the cover ${\mathcal Z}(\lambda)$. In particular $U_{i_m}$ will denote the cylinder $Z_{\lambda}(C_{i_m})$ associated with the vertex $v_{i_m} ^{\sigma}$ of $K$ for $0 \le m \le l$. Define the intersection ${\mathcal U} _{\sigma}$ to then be $U_{i_0} \cap\cdots\cap U_{i_l}$. Let $J_{\sigma} =
\{ j \in I - I_{\sigma} : U_j \cap {\mathcal U}_{\sigma} \neq \emptyset \}$. Define the set $V_{j,\sigma} = \{U_j \cap {\mathcal U}_{\sigma} : j \in J_{\sigma}\}$ and the family ${\mathcal V}_{\sigma}=(V_{j,\sigma})_{j \in J_{\sigma}}$.

We proceed to show that: 
\Claim ${\mathcal V}_{\sigma}$ is a cover for ${\mathcal U}_{\sigma}$. \EndClaim

\Proof Suppose on the contrary that ${\mathcal V}_{\sigma}$ is in fact \emph{not} a cover for
${\mathcal U}_{\sigma}$. Then there exists a point $P$ of ${\mathcal U}_{\sigma}$
such that $P \not\in U_i$ for any $i\in I-I_{\sigma}$. In particular, $G_P(\lambda) \le \theta(\sigma)$. However by \ref{thetasigma subgroup of G_P} we also have $\theta(\sigma) \le G_P(\lambda)$, and so $\theta(\sigma)= G_P(\lambda)$. Then because $r_{\theta(\sigma)} \le k-3$, we have $r_{G_P(\lambda)} \le k-3$. But, the minimum enveloping rank of $G_P(\lambda)$ is $\ge k-2$ as $r_M\ge k-2$
by hypothesis, providing a contradiction. Therefore, ${\mathcal V}_{\sigma}$ covers ${\mathcal U}_{\sigma}$ as claimed.
\EndProof

So ${\mathcal V}_{\sigma}$, which inherits the subspace topology, is in fact a cover of ${\mathcal U} _{\sigma}$, and so it follows from the definitions that the nerve of ${\mathcal V}_{\sigma}$ is simplicially isomorphic to the link of $\sigma$ in $K$. Note that two
different indices in $J_{\sigma}$ may define the
same set in ${\mathcal V}_{\sigma}$ but they will define different sets in
${\mathcal Z}(\lambda)$; this is why it is essential to define the nerve of
${\mathcal V}_{\sigma}$ using $J_{\sigma}$: so that the map from the vertex
set of the nerve of ${\mathcal V}_{\sigma}$ to the vertex set of the link of
$\sigma$ in $K$ is not only simplicial but bijective; that the
inverse of this map is simplicial is straightforward. 
To see this, suppose $v_j$ is a vertex in the nerve of ${\mathcal V}_{\sigma}$, then by
definition $U_j \cap {\mathcal U}_{\sigma} \neq \emptyset$, i.e.
$(U_{i_0} \cap \cdots \cap U_{i_l}) \cap U_j$ is non-empty. In
particular, $U_{i_0} \cap U_j$, $U_{i_1} \cap U_j, \dots ,
U_{i_l} \cap U_j$ are all non-empty, so that $\{ v_{i_0},
v_j \}, \{ v_{i_1}, v_j \}, \dots ,\{ v_{i_l}, v_j \}$
are all edges of $K$ ($v_j$ is distinct from the vertices of
$\sigma$), and $v_j$ is in the link of $\sigma$ in $K$. The
reverse inclusion is similar.  

Applying Borsuk's Nerve Theorem to ${\mathcal V}_{\sigma}$ in place of
${\mathcal Z}$, we see the underlying space of the nerve of
${\mathcal V}_{\sigma}$ is homotopy-equivalent to ${\mathcal U}_{\sigma}$. Since ${\mathcal U}_{\sigma}$ is a finite, non-empty intersection of convex open sets, it is
contractible. We conclude that the link in $K$ of every simplex of
minimum enveloping rank $m$ with $0 \le m \le k-3$ is contractible and non-empty.

We now show that the inclusion $|K|-|K_{(k-3)}| \rightarrow |K|$ is a
homotopy equivalence.

By local finiteness of the cover ${\mathcal Z}$ from which its nerve $|K|$ is defined, we may index the vertices of $K_{(k-3)}$, and therefore we may index the simplices of $K_{(k-3)}$ and partially order them in the following way: if $\sigma_i, \sigma_j$ are such that $\sigma_i$ is a proper face of $\sigma_j$, then $j<i$.

Define $F_n = \sigma_1 \cup \dots \cup \sigma_n$. We may regard $|K|-|K_{(k-3)}|$ as the topological direct limit of the subspaces $K_{F_n}=(|K|-|K_{(k-3)}|) \cup |F_n|$. Thus it suffices to show that the inclusion $K_{F_n} \hookrightarrow K_{F_{n+1}}$ 
is a homotopy equivalence. 


Now $K_{F_{n+1}} -  \sigma_{n+1} = K_{F_n}$ and as ${\rm lk}_{K}(\sigma_{n+1})$ is contractible by our work above, we may apply Proposition \ref{homoequiv} to get $K_{F_{n+1}} -  \sigma_{n+1} \cong K_{F_{n+1}}$. Hence the inclusion $K_{F_n} \rightarrow K_{F_{n+1}}$ is a
homotopy equivalence as required. \EndProof


\section{Group-Theoretic Preliminaries}
 \label{Group-Theoretic Stuff}

We will say that $W$ is a \emph{saturated} subset of the geometric realization $|K|$ of a simplicial complex $K$, if $W$ (endowed with the subspace topology) is a union of open
simplices of $|K|$ (endowed with the weak topology). 

Given a $\Gamma$-labeled complex $(K,(C_v)_v)$ and saturated subset $W \subseteq |K|$, we define the subgroup $\Theta(W)$ of $\Gamma$ to be the group $\la C_v : v \in \sigma, \sigma \subset W \ra$.

We now restate Group-Theoretic Conjecture \ref{two} from the Introduction which is necessary to prove Proposition \ref{meatandpotatoes general}, an essential ingredient in the proof of the Implication Theorem \ref{implication}. Let $H \vee K = \la H,K \ra$. 

\Conjecture \label{rankorama general} Suppose $H,K$ are two rank $h$ subgroups of a free group with $h \ge 3$. If the rank of $H \cap K$ is greater than or equal to $h$, then the rank of $H \vee K$ must be less than or equal to $h$.\EndConjecture

\Definition \label{local rank} We say a group $\Gamma$ has \emph{local rank $\le k$} where $k$ is a positive integer, if every finitely generated subgroup of $\Gamma$ is
contained in a subgroup of $\Gamma$ which has rank less than or equal
to $k$. The local rank of $\Gamma$ is the smallest $k$ with this
property. If there does not exist such a $k$ then we define the local
rank of $\Gamma$ to be $\infty$. Note that if $\Gamma$ is
finitely generated, its local rank is simply its rank. \EndDefinition

\Proposition \label{meatandpotatoes general} Assume Conjecture \ref{rankorama general}. Let $k,r \in {\bf Z}^+$ with
$k>r \ge 3$ and $k \ge 5$. Suppose $\Gamma$ is a $k$-free group,
$(K, (C_v)_v)$ a $\Gamma$-labeled complex, and $W$ a saturated, connected subset of $|K|$ such that ${\rm rk} \; \Theta(\sigma) =
r$ for all $\sigma \subset W$. Assume additionally that either \Alternatives

\item \label{meatandpotatoes1 general} there exists a positive integer $n$
  such that for all open simplices $\sigma$ in $W$, the dimension of $\sigma$ is $n$ or $n-1$, or

\item \label{meatandpotatoes2 general} $r=k-2$ and $\sigma \in
|K^{(k-1)}|-|K_{(k-3)}|$ for all $\sigma \in W$. \EndAlternatives 
Then the local rank of $\Theta(W)$ is at most $r$. \EndProposition

\Proof By definition, we are required to show that every finitely generated subgroup of
$\Theta(W)$ is contained in a finitely generated subgroup of $\Theta(W)$ which has rank less than or equal to $r$. So suppose that $E \le \Theta(W)$ is a finitely
generated subgroup of $\Theta(W)$. Then $E\le \Theta(V_0)$ for some
saturated subset $V_0$ of $W$ that contains finitely many open
simplices. Because $W$ is connected and $V_0$ contains only finitely
many open simplices, there is a smallest
connected subset $V$ of $W$ that is a union of finitely many open
simplices such that $V_0 \subseteq V$; clearly $E \le \Theta(V)$ and
$V$ is finitely generated. We will show by induction on the number of
simplices in $V$ that $\Theta(V)$ has rank at most $r$.

Proceeding as in \cite[Proposition 4.4]{cs}, by connectedness we
may list the (finitely many) open simplices of $V$ in the
following way: $\sigma_0,\dots,\sigma_m$, ($m \ge 0$ since $V$ is
non-empty) where for any $i$ with $0 \le i \le m$, there is an
index $l$ with $0 \le l < i$ such that $\sigma_l$ is a proper face
of $\sigma_i$ or $\sigma_i$ is a proper face of $\sigma_l$. Define
the saturated subset $V_i = \sigma_0 \cup \dots \cup \sigma_i$ for
$0 \le i \le m$; by induction on $i$, we will show ${\rm rk}
\;\Theta(V_i) \le  r$. The base case is straightforward as
$\Theta(V_0)= \Theta(\sigma_0)$ and $\sigma_0$ is an open
simplex of $W$, and so has rank $r$ by hypothesis. For the induction step assume ${\rm rk}\;\Theta(V_{i-1})=r$; we want to show that ${\rm rk}\;\Theta(V_i)=r$. By how we have arranged the list of simplices in $V$, there is an index
$l$ with $0 \le l < i$ such that $\sigma_l$ is a proper face of
$\sigma_i$ or $\sigma_i$ is a proper face of $\sigma_l$.

Case (i): First consider the case when $\sigma_i$ is a proper face
of $\sigma_l$. Then $\Theta(V_i) = \Theta(V_{i-1})$ as $\sigma_i <
\sigma_l \in V_{i-1}$. By our induction assumption, ${\rm rk}
\;\Theta(V_{i-1}) \le r$, and so ${\rm rk} \;\Theta(V_i) \le r$ as
required.

Case (ii): Next, consider the case when $\sigma_l$ is a proper
face of $\sigma_i$. Let $P=\Theta(V_{i-1})$, $Q=\Theta(\sigma_i)$
and $R=\Theta(\sigma_l)$. Then ${\rm rk} \;P \le k-2$ by the induction hypothesis and ${\rm rk}
\;Q = {\rm rk} \;R = k-2$ by assumption. We want to show $\Theta
(V_i) = P \vee Q$ has rank less than or equal to $r$.

Subcase (i): Assume first that property (\ref{meatandpotatoes1 general}) holds.
Then since $\sigma_l$ is a proper face of
$\sigma_i$, we must have ${\rm dim}\sigma_i=n$ and ${\rm dim}\sigma_l
= {\rm dim}\sigma_i-1=n-1$. Let $v$ denote the vertex of $\sigma_i$
such that ${\rm span}\{ \sigma_l, v \} = \sigma_i$ and let $C=C_v$. Then $Q = R \vee C$, and $P \vee C = P \vee Q$. So we
proceed to show that ${\rm rk} (P \vee C) \le r$.

By way of contradiction, assume ${\rm rk}(P \vee C) > r$.
Then since $C$ is infinite cyclic, $P \vee C$ has rank at
most ${\rm rk} \;P +1= r+1$ and so $P \vee C$ has rank
exactly $r+1$. As $\Gamma$ is $k$-free and $r<k$ (and hence $r+1
\le k$), it follows that $P \vee C$ is free as a subgroup
of $\Gamma$ \emph{and} in particular is the free product of the subgroups
$P$ and $C$ (\cite[Lemma 4.3]{cs}). But, since $R \le P$, in
particular $Q=R \vee C$ is the free product of $R$ and
$C$, and so has rank equal to ${\rm rk} \;R+1 = r+1$, which is a
contradiction as the rank of $Q$ is exactly $r$. We conclude that
$P \vee C $ has rank $\le r$ as required for this subcase.

Subcase (ii): Next we assume property (\ref{meatandpotatoes2 general}). Then
$r= {\rm rk} \;Q={\rm rk} \;R = k-2$ (${\rm rk} \;P \le k-2$ by induction assumption). As
$r= {\rm
  rk} \;Q={\rm rk} \;R = k-2$, both ${\rm dim} \sigma_l$ and
${\rm dim} \sigma_i$ are at least $k-3$. Also $\sigma_l, \sigma_i \in
K^{(k-1)}$, so both ${\rm dim} \sigma_l$ and ${\rm dim} \sigma_i$ are $\le k-1$. Finally,
since our Case (ii)-assumption is that $\sigma_l$ is a proper face of
$\sigma_i$, possible pairs $({\rm dim} \sigma_l , {\rm dim} \sigma_i)$
are $(k-3,k-2)$, $(k-2,k-1)$, and $(k-3,k-1)$. Let $C \le \Gamma$ denote the
subgroup of $Q$ such that $Q = R \vee C$; then $P \vee Q = P \vee C$ since $R \le P$.

\Number \label{norank2foryou general} First, we look at the rank of
$P$. A priori we know that ${\rm rk}\;P \ge 2$
(i.e. $P$ cannot be cyclic and is non-trivial) since $P$ contains
the rank-$(k-2)$ subgroup $R$. 

In particular, as $R= \Theta(\sigma_l)$ is a subgroup of $P$, and as
$\sigma_l$ is an element of $K^{(k-1)}-K_{(k-3)}$, we know that the minimum enveloping rank of $R$ is strictly greater than $k-3$. Along
with our induction assumption that ${\rm rk}\;P \le k-2$, we conclude
the rank of $P$ is exactly $k-2$. (Note that for this reason in the case when $k=4$,
it is enough only to say in (\ref{meatandpotatoes2 general}) that $r=2$,
since what is required for the rest of the argument is that $P$ have
rank exactly $2=k-2$ in this case, an immediate consequence of $P$
containing the rank $2$ subgroup $R$. Specifically, in the $k=4$ case, we see that a group containing a rank
two subgroup certainly cannot have rank one; whereas in cases for $k \ge 5$, one
observes that a group that contains a rank three (or more) subgroup \emph{can} have
rank two or more, and so that $r_{\theta(\sigma)} \ge k-2$ is required
in the statement of (\ref{meatandpotatoes2 general})). Next observe that we must have ${\rm rk}\;C = 1$ or $2$ as demonstrated by
the possible pairs $({\rm dim} \sigma_l , {\rm dim} \sigma_i)$ above.
All together, this gives that ${\rm rk} (P \vee C) \le k$ and so $P \vee C$ is free as
a subgroup of $\Gamma$. \EndNumber


\Number \label{norank2foryou2 general} Next, notice that because $Q
\le P \vee Q$ and $Q$ has minimum enveloping rank $\ge k-2$,
$P \vee Q$ cannot have rank less than $k-2$. Along with the bound ${\rm rk} (P \vee C) \le k$ of \ref{norank2foryou general}, we conclude there are only three possibilities for
the rank of the group $P \vee Q (=P \vee C)$: these are $k, k-1, {\rm and}\; k-2$.\EndNumber

\Number \label{norank2foryou3 general} As we have $R \le P$, $R \le Q$, and $R \le P \cap Q$, then for the same reason as outlined in \ref{norank2foryou general} with $P \cap Q$ taking the
place of $P$, we conclude ${\rm rk} ( P \cap Q )
>k-3$. Therefore, we may apply Conjecture \ref{rankorama general} which gives
that ${\rm rk} (P \vee Q) \le k-2$, and so must be equal to
$k-2$ by \ref{norank2foryou2 general}, completing this final subcase
and proving the proposition. \EndNumber \EndProof



 
\section {Theorem and general bound on $r_M$} \label{Theorem and Corollary}

We now restate formally and prove the implicative statement of \ref{implication} given in the Introduction. For the proof we require a few basic definitions about graphs.

\Definitions \label{graph and tree} We say that ${\mathcal G}$ is a \emph{graph} if ${\mathcal G}$ is at most a one-dimensional simplicial complex (and so ${\mathcal G}$ has no loops or multiple edges). A \emph{tree} $T$ is a connected graph with no cycles; i.e. $T$ is a graph which is simply connected. Further, if $X_i$ and $X_j$ are disjoint, saturated subsets of a simplicial complex $|K|$, we will make use of the concept of an abstract bipartite graph ${\mathcal G}={\mathcal G}(X_i,X_j)$ constructed in the following way. Let ${\mathcal W}_i$, ${\mathcal W}_j$ be the sets of connected components of $X_i$ and $X_j$ respectively. Then the vertices of ${\mathcal G}$ are the elements of ${\mathcal W}_i \cup {\mathcal W}_j$, and a pair $\{ v_{W_i}, v_{W_j} \}$ is an edge if there exist simplices $\sigma \in W_i$ and $\tau \in W_j$ for which $\sigma \le \tau$ or $\tau \le \sigma$. Finally, we say that the simplicial action of a group $\Gamma$ on a graph ${\mathcal G}$ is without inversions if for every $\gamma \in \Gamma$ that stabilizes an edge $e= \{ v_1,v_2 \} \in {\mathcal G}$, we have $\gamma\cdot v_1 = v_1$ and $\gamma\cdot v_2 = v_2$. \EndDefinitions

The following two lemmas taken directly from \cite{cs} will provide the contradiction necessary to prove Theorem \ref{implication}:

\Lemma \label{Shalen's 2} Suppose that $K$ is a simplicial complex and that $X_i$ and $X_j$ are saturated subsets of $|K|$. Then $|{\mathcal G}(X_i, X_j)|$ is a homotopy-retract of the saturated subset $X_i \cup X_j$ of $|K|$. \EndLemma
\Proof This is \cite[Lemma 5.12]{cs}. \EndProof

\Lemma \label{Shalen's 1} Let $M$ be a closed, orientable, aspherical $3$-manifold. Then $\pi_1(M)$ does not admit a simplicial action without inversions on a tree $T$ with the property that the stabilizer in $\pi_1(M)$ of every vertex of $T$ is a locally free subgroup of $\pi_1(M)$.\EndLemma
\Proof This is \cite[Lemma 5.13]{cs} \EndProof


Finally, we will appeal to the property stated in the next remark in the proof of Theorem \ref{rM is less than k-3 again}.

\Remark \label{stabnorm} Suppose a group $\Gamma$ admits a labeling-compatible action on a $\Gamma$-labeled complex $(K, (C_v))_v$, as is defined in \ref{labeling-compatible}. If $W$ is a saturated subset of $|K|$ and $\gamma$ is any element of $\Gamma$, it follows that
$\Theta (\gamma\cdot W) = \gamma \Theta(W) \gamma^{-1}$. (Since by the definitions, $\Theta (\gamma\cdot W) = \la C_v : v \in
\gamma \cdot W \ra = \la C_{\gamma \cdot v} : v \in W \ra = \la
\gamma C_v \gamma^{-1} : v \in W \ra = \gamma \la C_v : v \in W
\ra \gamma^{-1} = \gamma \Theta(W) \gamma ^{-1}$). So if an
element $\gamma$ of $\Gamma$ is invariant on $W$, then it is in the
normalizer of $\Theta(W)$. More generally, the stabilizer in
$\Gamma$ of $W$ is a subgroup of the normalizer of $\Theta(W)$.
\EndRemark

The following theorem is the reformulated Implication Theorem \ref{implication} of the Introduction.

\Theorem \label {rM is less than k-3 again} Suppose $M$ is a closed, orientable, hyperbolic 3-manifold such that $\pi_1 (M)$ is $k$-free with $k\ge 5$. Then if one assumes the Conjecture of \ref{two} with $m=k-2$, setting $\lambda = \log{(2k-1)}$ we have $r_M \le k-3$.\EndTheorem

\Proof We have $M = {\bf H}^3/\Gamma$, where $\Gamma \le {\rm Isom}_+
({\bf H}^3)$ is discrete, compact, and torsion-free.

We will assume that $r_M \ge k-2$ and proceed by way of
contradiction. Equivalently, suppose that for all points $P$ in
${\bf H}^3$, the minimum enveloping rank of $G_P$ is $\ge k-2$. Then in particular, ${\bf
H}^3 = \bigcup_ {C_1^P,\dots,C_{k-2}^P \in C_{\log{(2k-1)}}(\Gamma), P\in
{\bf H}^3}
Z_{\log{(2k-1)}} (C_1^P) \cap \cdots \cap Z_{\log{(2k-1)}} (C_{k-2}^P)$ as described in \ref{there is a point}. Without loss of generality we write 

${\bf H}^3 = \bigcup_ {C_1,\dots,C_{k-2} \in C_{\log{(2k-1)}}(\Gamma)}
Z_{\log{(2k-1)}} (C_1) \cap \cdots \cap Z_{\log{(2k-1)}} (C_{k-2})$, and define the family

${\mathcal Z}=
(Z_{\log{(2k-1)}} (C_i))_{C_i \in C_{\log{(2k-1)}}(\Gamma), 1\le i \le k-2}$.

We have that ${\mathcal Z}$ is an
open cover of ${\bf H}^3$ which satisfies the hypothesis of Lemma \ref{rank argument
general}. Then if $K$ denotes the nerve of ${\mathcal Z}$, the result gives that $|K|-|K_{(k-3)}| \cong {\bf H}^3$. Since the
inclusion $|K^{(n)}|-|K_{(k-3)}| \rightarrow |K|-|K_{(k-3)}|$
induces isomorphisms on $\pi_0$ and $\pi_1$ for $n\ge k-1$ (see \cite[Lemma 5.6]{cs}), 
it follows that $|K^{(k-1)}|-|K_{(k-3)}|$ is connected and simply
connected.

Let $\sigma$ be an open simplex in $|K^{(k-1)}|-|K_{(k-3)}|$. Applying
Lemma \ref{the first one} with $n={\rm dim}(\sigma) + 1$ (i.e. $n$ is the
number of vertices of $\sigma$ and therefore the number of associated maximal cyclic subgroups of $\Gamma$ whose associated cylinders have nonempty intersection, as is determined by the nerve), 
we have that the rank of $\Theta(\sigma)$ is less than or equal to $k-1$. Now
since $\sigma$ is in $|K^{(k-1)}|-|K_{(k-3)}|$, by definition the minimum enveloping rank
of $\Theta(\sigma)$ is at least $k-2$. In particular, the rank of $\Theta(\sigma)$ is at least $k-2$ by \ref{about mer}.

\Number \label{disjoint union} All together, we conclude that for
any open simplex $\sigma$ in $|K^{(k-1)}|-|K_{(k-3)}|$, the rank of
$\Theta(\sigma)$ is $k-2$ or $k-1$. So, we may write
$|K^{(k-1)}|-|K_{(k-3)}|$ as a disjoint union of the saturated subsets $X_{k-2}$ and
$X_{k-1}$, where $X_i$ for $i=k-2, k-1$ is the union of all open simplices $\sigma$ of
$K^{(k-1)}$ for which $\Theta(\sigma)$ has rank $i$. \EndNumber

We claim the following: 
\Claim \label{localrank} For $i \in \{ k-2,k-1 \}$ and for any component $W$ of $X_i$, the local rank of $\Theta(W)$ is at most $i$. \EndClaim

\Proof First we consider the case when $i=k-2$. Then $W$ is a
component of $X_{k-2}$, and for any open simplex $\sigma$ of $X_{k-2}$, ${\rm rk}\;\Theta (\sigma)$ is exactly $k-2$. 
Taking $r=k-2$ in Proposition \ref{meatandpotatoes general},
specifically in item (\ref{meatandpotatoes2 general}), we have that the local rank of $\Theta(W)$ is
at most $k-2$. The proof in the case of (\ref{meatandpotatoes2 general})
shows that the local rank of $\Theta(W)$ is \emph{exactly} $k-2$. 

Suppose next that $i=k-1$. Then $W$ is a component of $X_{k-1}$, and so for each open simplex $\sigma$ of $X_{k-1}$, we have $\rank \Theta (\sigma)$ is exactly $k-1$. If $d$ denotes the dimension of $\sigma$, then the
subgroup $\Theta (\sigma)$ is generated by $d+1$ cyclic groups
which are elements of $C_{\log{(2k-1)}}(\Gamma)$. Hence ${\rm rk}\;\Theta(\sigma) \le d+1$ and in particular $d \ge {\rm rk}\;\Theta (\sigma) -1$. As  ${\rm rk} \; \Theta(\sigma) =k-1$, we have $d \ge k-2$. But because $\sigma$ is a simplex contained in $K^{(k-1)}$, $d$ is less than or equal to $k-1$, and so we must have $d=k-2$ or $k-1$. Letting $r=k-1$ and $n=k-1$ in item (\ref{meatandpotatoes1 general}) of Proposition \ref{meatandpotatoes general}, we satisfy the hypotheses and the
conclusion gives that $\Theta(W)$ has local rank at most $r=k-1$ as desired. \EndProof

Next, we claim: \Claim \label{localrankcorollary} The local rank
of $\Theta(W)$ is exactly $k-2$ or $k-1$. \EndClaim

\Proof Let $l_W$ be the \emph{local} rank
of $\Theta(W)$. Our previous claim shows that $l_W \le k-1$. If in
fact $l_W \le k-3$, then by definition any finitely generated
subgroup of $\Theta(W)$ is contained in a finitely generated
subgroup of rank less than or equal to $k-3$. As $\Theta(\sigma) \le \Theta(W)$, this says that $\Theta(\sigma)$ is contained in a subgroup of rank less than or equal to $k-3$ and so the minimum enveloping rank of $\Theta(\sigma)$ would be $\le k-3$ in this situation. However, given an open
simplex $\sigma$ in $W$, in particular $\sigma$ is a simplex of
$|K^{(k-1)}|-|K_{(k-3)}|$ and so $\Theta(\sigma)$ has minimum enveloping rank $\ge k-2$, providing a contradiction.
Therefore, $l_W$ is $k-2$ or $k-1$. \EndProof

\Claim \label{normalrank}(The analogue of \cite[Claim 5.13.2]{cs})
If $W$ is a component of $X_{k-2}$ or $X_{k-1}$, the normalizer of
$\Theta(W)$ in $\Gamma$ has local rank at most $k-1$. \EndClaim

\Proof As a subgroup of $\Gamma$, the normalizer of $\Theta(W)$ is
$k$-free. Clearly $\Theta(W)$ is a normal subgroup of its
normalizer, and since by the result of \ref{localrankcorollary} we have $l_W =
k-2$ or $k-1$ which are strictly less than $k$, it follows by \cite[Proposition 4.5]{cs} 
that the normalizer of $\Theta(W)$ has local rank at most $l_W$. \EndProof

Set $T= {\mathcal G} (X_{k-2}, X_{k-1})$ (see Definitions \ref{graph and tree}). By Lemma \ref{Shalen's 2}, $T$ is a
homotopy-retract of $X_{k-2} \dot\cup X_{k-1}$, which is equal to $|K^{(k-1)}|-|K_{(k-3)}|$ by \ref{disjoint union}.
Since $|K^{(k-1)}|-|K_{(k-3)}|$ is connected and simply connected, $T$ is a tree.

By Definition \ref{labeling-compatible} of the $\Gamma$-labeling compatible action of
$\Gamma$ on $K$, we see that for any $\gamma \in \Gamma$ and
$\sigma$ in $K^{(k-1)}$, $\Theta(\sigma)$ and
$\Theta(\gamma\cdot\sigma)$ are conjugates in $\Gamma$ (see Remark \ref{stabnorm}), and so have equal rank. Consequently, $X_{k-2}$ and $X_{k-1}$ are invariant under
the action of $\Gamma$. Note that if $w$ is a vertex of $T$, the
stabilizer $\Gamma_w$ of $w$ in $\Gamma$ is really the stabilizer
of the associated component $W$ in $X_{k-2}$ or $X_{k-1}$, and so by Remark \ref{stabnorm}, $\Gamma_w \le {\rm normalizer}\;\Theta(W)$.

\Number \label{freeatlast} By our work above in \ref{normalrank},
the local rank of ${\rm normalizer}\;\Theta(W)$ is at most $k-1$,
and given that it contains $\Gamma_w$ as a subgroup, $\Gamma_w$
must also have local rank at most $k-1$, and, in particular, is
locally free being a subgroup of $\Gamma$ which is $k$-free.
\EndNumber

Therefore we've constructed an induced action by $\Gamma$ on the
tree $T$ \emph{without} inversions. Since the stabilizer of any vertex of
$T$ is locally free as a subgroup of $\Gamma$ by \ref{freeatlast},
our construction admits a contradiction to Lemma \ref{Shalen's 1}.
\EndProof

The following Propositions and Definitions will be used to explain the geometry of the cases when $r_M(\lambda)=0$ and $1$, and in particular will be used when $\pi_1(M)$ is $5$-free and $\lambda = \log9$ in Corollary \ref{cor for alts for hyp sp}.

\Proposition \label{rank zero} Suppose $\lambda >0$ and $M = {\bf
H}^3 / {\Gamma}$ is a closed, orientable hyperbolic 3-manifold
with $\Gamma$ discrete and purely loxodromic. If $r_M = 0$, then M
contains an embedded ball of radius $\lambda / 2$. \EndProposition

\Proof  As $r_M=0$, ${\rm rk}\;G_P \ge 0$ for all $P \in {\bf
H}^3$, and in particular, the choice of $r_M$ means there is a
point $P_0 \in {\bf H}^3$ with ${\rm rk}\;G_{P_0} = 0$. Then $P_0 \not\in
Z_{\lambda}(C)$ for any $C \in C_{\lambda}(\Gamma)$, and so
$d(P_0, \gamma \cdot P_0) \ge \lambda$ for all $\gamma \in \Gamma
- \{ 1 \}$, and more generally, ${\bf H}^3 \not = \bigcup_{C \in
C_{\lambda}(\Gamma)} Z_{\lambda}(C)$. If $B_{P_0}(\lambda / 2)$
denotes the hyperbolic open ball of radius $\lambda / 2$ with
center $P_0$, in particular this says that the injectivity radius
of $B_{P_0}(\lambda / 2)$ in $M$ is $\lambda / 2$; namely
$B_{P_0}(\lambda) \cap \gamma \cdot B_{P_0}(\lambda) =
\emptyset$. To see this, consider a point $P'$ in
$B_{P_0}(\lambda)$. If in fact it was true that $\gamma(P')$ is
also in $B_{P_0}(\lambda)$, it would then follow that
$d(P_0,\gamma \cdot P_0) \le d(P_0,\gamma \cdot P') +
d(\gamma \cdot P',\gamma \cdot P_0) < \lambda /2 + \lambda /2 = \lambda$,
giving a contradiction. Therefore if $q:{\bf H}^3 \rightarrow M$
is the projection map, $q|B:B\rightarrow M$ is injective and the
conclusion follows. \EndProof

\Definitions \label{mathfrak X} Let ${\mathfrak X}_M$ be the set of
points $P$ in $M$ such that if $l_P$ denotes the length of the
shortest, homotopically non-trivial loop based at $P$, then there
is a maximal cyclic subgroup $D_P$ of $\pi (M,P)$ such that for
every homotopically non-trivial loop $c$ based at $P$ of length
$l_P$, we have $[c] \in D_P$. Note that the loop $c$ of length $l_P$ may represent a proper power of a generator of $D_P$. Let ${\mathfrak s}_M(P)$ be the smallest length of any loop $c$ based at $P$ such that $[c] \not \in D_P$.\EndDefinitions

\Proposition \label{rank one} Suppose $\lambda >0$ and $M = {\bf
H}^3 / {\Gamma}$ is a closed, orientable hyperbolic $3$-manifold
with $\Gamma$ discrete and torsion-free. If $r_M = 1$, then there
exists a point $P^{*} \in {\bf H}^3$ with $P^{*} \in {\mathfrak X}_M$ and ${\mathfrak s}_m(P^{*}) =\lambda$.
\EndProposition

\Proof As $r_M = 1$, the definition of $r_M$ gives that ${\rm
rk}\;G_P \ge 1$ for all $P \in {\bf H}^3$ (more generally that
${\bf H}^3 = \bigcup_{C \in C_{\lambda}(\gamma)} Z_{\lambda}(C)$) \emph{and} that there is a point $P_0\in {\bf H}^3$ with ${\rm
rk}\;G_{P_0} = 1$. Hence $P_0 \in Z_{\lambda}(C_0)$ for some $C_0
\in C_\lambda(\Gamma)$ and $P_0 \not \in Z_{\lambda}(C)$ for any
other $C \in C_{\lambda}(\Gamma)-C_0$, namely $G_{P_0} = \la C_0
\ra$. Set $Z_0= Z_{\lambda}(C_0)$ and $Y= \bigcup_{C \in
C_\lambda(\Gamma)- C_0} Z_{\lambda}(C)$. Then ${\bf H}^3 = Y \cup
Z_0$. Since ${\bf H}^3$ is connected, $(Z_{\lambda}(C))_{C \in C_\lambda(\Gamma)}$ is an open cover, and
$\Gamma$ is discrete, we must have the intersection $Y \cap Z_0$ is nonempty and open. Notice $P_0 \not \in Y$ means $Z_0 \not \subseteq Y$. As $Z_0$ is connected, we conclude that the frontier of the set $Y \cap Z_0$ relative to $Z_0$ is nonempty; let $F$ denote this set. Let us choose a point $P^{*}$ in $F$. In particular, this says
that (i) $P^{*} \in Z_0$ \emph{and} (ii) $P^{*}$ is in the frontier of $Y$
(relative to ${\bf H}^3$). (In concluding (ii), recall that the
collection of cylinders in $Y$ comprises a locally finite
collection because $\Gamma$ is \emph{discrete}, and so $P^{*}$, a limit
point of $Y$, does not belong to this open collection). If
$\gamma_0$ is a generator for $C_0$, (i) implies that $d(P^{*},
\gamma_0^m \cdot P^{*}) < \lambda$ for some integer $m \ge 1$. By
(ii), we know that $d(P^{*}, \gamma_1 \cdot P^{*}) = \lambda$ for some
$\gamma_1 \in \Gamma - \gamma_0$ \emph{and} that $d(P^{*}, \gamma \cdot P^{*})
\ge \lambda$ for all $\gamma \in \Gamma - C_0$. Using the base
point $P^{*} \in {\bf H} ^3$ to identify $\pi (M, q(P^{*}))$ with
$\Gamma$, we have that $\gamma_0^m$ is represented by a loop of
length less than $\lambda$ based at $q(P^{*})$, and any other homotopically non-trivial loop of length less than $\lambda$ based at $q(P^{*})$ is identified with an element of $C_0$.
Therefore, we have shown the existence of a point $P^{*} \in
{\mathfrak X}_M$ for which the smallest length of any loop represented by $[c]$ in
$M$ based at $P^{*}$ with the property that $[c]$ is not in $D_{P^{*}}$, is exactly
$\lambda$. \EndProof


\section{Matrices and Theorem for the case k=5}\label{Matrices}
We will now restate some of Kent's constuction and results regarding joins and intersections of free groups; in particular, we incorporate the background (\ref{graphs background} and \ref{chi}) as discussed in \cite{Kent} which is needed to apply \cite[Lemma 7]{Kent} and \cite[Lemma 8]{Kent} to prove Proposition \ref{Kent's} that follows. Subsequently, Conjecture \ref{two} for rank-$3$ subgroups $H$ and $K$ is affirmed in Corollary \ref{Corollary for 5}.

\Number \label {graphs background} Given a free group $F$ free on the set $\{ a,b \}$, we associate with $F$ the wedge $W$ of two circles based at the wedge point, and we orient the (two) edges of $W$. Then for any subgroup $H$ of $F$ there is a unique choice of basepoint $*$ in the covering space $\widetilde{W_H}$ such that $\pi_1(\widetilde{W_H}, *)$ is exactly $H$.  Then $\Gamma_H$ will denote the smallest subgraph containing $*$ of $\widetilde{W_H}$ that carries $H$. By this construction, $\Gamma_H$ inherits a natural oriented labeling, i.e. each edge is labeled with one of $\{ a,b \}$ and initial and terminal vertices (not necessarily distinct) are determined by the orientation. Hence ${\rm rk}\pi_1(\Gamma_H) = {\rm rk}\;H$. Also by this construction, any vertex of $\Gamma_H$ is at most $4$-valent. Define a $3$- or more valent vertex of $\Gamma_H$ to be a \emph {branch vertex}. We will from here on assume that all graphs in our duscussion are normalized so that all branch vertices are $3$-valent (see the beginning of \cite[Section 3]{Kent} for this explanation). \EndNumber

\Number \label{chi} If $\Gamma$ is a graph, let $b(\Gamma)$ denote the number of branch vertices in $\Gamma$. Note that if $\Gamma$ is $3$-regular, i.e. all branch vertices are $3$-valent, then we have $\overline{\chi}(\Gamma)={\rm rk}(\pi_1(\Gamma))-1=b(\Gamma)/2$. By \ref{graphs background} this says that if $H,K$ are subgroups of $F$, then ${\rm rk} (\pi_1(\Gamma_{H \vee K}))-1= {\rm rk} (H \vee K) -1$. If $V(\Gamma_H)$ and $V(\Gamma_K)$ denote the vertex sets of $\Gamma_H$ and $\Gamma_K$ respectively, then we can define the graph ${\mathcal G}_{H \cap K}$ whose vertex set is the product $V(\Gamma_H) \times V(\Gamma_K)$ and for which $\{ (a,b),(c,d) \}$ is an edge if there are edges $e_1=\{ a,c \}$ in $\Gamma_H$ and $e_2=\{ b,d \}$ in $\Gamma_K$ for which $e_1$ and $e_2$ have the same label, and $e_1$ is oriented from $a$ to $c$ and $e_2$ is oriented from $b$ to $d$. The graph ${\mathcal G}_{H \cap K}$ is the pullback of the maps $\Gamma_H \rightarrow W$ and $\Gamma_K \rightarrow W$ in the category of oriented graphs, and $\Gamma_{H \cap K}$ is a subgraph of ${\mathcal G}_{H \cap K}$ that carries the fundamental group. We then have the projections  $\Pi_H: {\mathcal G}_{H \cap K} \rightarrow \Gamma_H$ and $\Pi_K: {\mathcal G}_{H \cap K} \rightarrow \Gamma_K$. Let the graph ${\mathcal T}$ denote the topological pushout of the maps $\Gamma_{H \cap K} \rightarrow \Gamma_H$ and $\Gamma_{H \cap K} \rightarrow \Gamma_K$ in the category of not properly labeled oriented graphs. Hence the graph ${\mathcal T}$ is defined as the quotient of the disjoint union $\Gamma_H \cup \Gamma_K$ modulo $\sim_{\mathcal R}$ where $x \in H$ is equivalent to $y \in K$ if $x \in \Pi_H(\Pi_K^{-1}(y))$ or $y \in \Pi_K(\Pi_H^{-1}(x)) \Pi_H$. Now since the map ${\mathcal T} \rightarrow \Gamma_{H \vee K}$ factors into a series of folds (which is surjective at the level of $\pi_1$), it follows that $\chi({\mathcal T}) \le \chi(\Gamma_{H \vee K})$. Equivalently $\overline{\chi}(\Gamma_{H \vee K}) \le \overline{\chi}({\mathcal T})$. 
\EndNumber

\Number \label{setup Kent's} 
As outlined in Kent's Section $3.2$ \cite{Kent}, we consider the $(2h-2)\times(2k-2)$ matrix, $M= (f(x_i,y_j))$, where $f:X\times Y\rightarrow \{ 0,1 \}$ is the function defined on the sets $X= \{x_1,\dots,x_{2h-2} \}$ and $Y=\{y_1,\dots,y_{2k-2} \}$ of branch vertices of $\Gamma_H$ and $\Gamma_K$, respectively, by letting $f(x_i,y_j)=1$ if $(x_i,y_j)$ is a branch vertex of $\Gamma_{H \cap K}$ and $0$ otherwise. 
Then the number of ones in $M$ is equal to the number of valence-$3$ vertices in $\Gamma_{H \cap K}$; i.e. $b(\Gamma_{H\cap K}) = \Sigma_{i,j} f(x_i,y_j)$. From \ref{chi}, we have $\overline{\chi}(\Gamma_{H \cap  K})=b(\Gamma_{H\cap K})/2$. By \cite[Lemma 8]{Kent}, after permuting rows and columns of $M$, we may write $M$ in the block form:
$(M_1, \dots, M_l, O_{(p \times q)})$ where $O_{(p \times q)}$ is the $p \times q$ zero matrix, every row and every column of each of the $M_i$ has a $1$. Here, every block $M_i$ of $M$ represents a unique equivalence class of $\sim_{\mathcal R}$ with representatives in $\Gamma_H$ and $\Gamma_K$; all-zero rows represent the $\le p$ equivalence class(es) of $\Gamma_H$ which do not have representatives in $\Gamma_K$; and all-zero column(s) represent the $\le q$ equivalence class(es) of $\Gamma_K$ which do not have representatives in $\Gamma_H$.
\EndNumber

\Proposition \label {Kent's} ${\rm rk} (H \vee K) \le 1+ \frac{1}{2} (l+p+q)$.\EndProposition
\Proof
Note that $2h-2 \ge l+p$ and $2k-2 \ge l+q$ as $2h-2$ is $\# \{ {\rm rows\;of} \; M \}$ and $2k-2$ is $\# \{ {\rm columns\;of} \; M \}$ where $M$ is the block matrix of \ref{setup Kent's}, and so $l$ is bounded above by the positive integer ${\rm min}((2h-2)-p,(2k-2)-q)$. We have ${\rm rk} (H \cup K) - 1 \le \overline{\chi}({\mathcal T}) \le \frac{1}{2} (l+p+q)$ by combining \ref{chi} and \ref{graphs background} along with \cite[Lemma 7]{Kent} for the first inequality and \cite[Lemma 8]{Kent} for the second. In particular, ${\rm rk} (H \vee K) \le 1+ \frac{1}{2} (l+p+q)$ as required.
\EndProof

\Corollary \label{Corollary for 5} Suppose $h=k=3$ and ${\rm rk} (H \cap K) \ge 3$. Then ${\rm rk} (H \vee K) \le 3$.
\EndCorollary

\Proof  As $h=k=3$, we consider the $4\times4$ block matrix $M$ of \ref{setup Kent's} where each row and each column of each of the $M_i$ has a $1$. So the number of ones, which is the number of valence-$3$ vertices in $\Gamma_{H\cap K}$, is $\ge4$. \Number \label{bound on l for 4} Note that $l$ is bounded above by ${\rm min}(4-p,4-q)$ by the proof of Proposition \ref{Kent's}, and so in particular, $l+p \le 4$ and $l+q \le 4$. Thus we may rewrite Proposition \ref{Kent's} to read ${\rm rk} (H \vee K) \le 1+{\rm min}(2+\frac{p}{2}, 2+\frac{q}{2})$. Using this formula, we see that ${\rm rk} (H \vee K) \le 3$ unless $p$ and $q$ are $\ge2$, and so we need only consider the following cases: \EndNumber

\emph {Case $p=4$ or $q=4$:} This case is an impossibility, as this would imply $M=O_{4 \times 4}$, and hence the number of branch vertices of $\Gamma_{H \cap K}$ is zero, a contradiction, and so we must have $p,q\le3$.

\emph {Case $p=3$ or $q=3$:} Suppose first that $p=3$. Then \ref{bound on l for 4} says that $l \le1$, implying that $l=1$ and the top row of $M$ has $4$ ones, and so $q=0$. In this case, the inequality of \ref{bound on l for 4} gives ${\rm rk} (H \vee K) \le 1+{\rm min}(2+\frac{3}{2}, 2)$, and so ${\rm rk} (H \vee K)\le 3$. When $q=3$ the argument is symmetric, and so the conclusion is satisfied.

\emph {Case $p=2$ or $q=2$:} By symmetry assume $p=2$. This gives $l \le 2$ by the bound on $l$ of \ref{bound on l for 4}. Now if $l\le1$, then $q$ must be equal to $3$ to satisfy the requirement on the number of ones in $M$ (i.e. the valence-$3$ vertices in $\Gamma_{H \cap K}$), which is the previous case. Next, if $l=2$, then as $l \le {\rm min}(4-p,4-q)$, we have $q\le2$. First, if $p=q=2$, then the requirement that the number of ones in $M$ is $\ge 4$ fails as the values of $p,q,$ and $l$ would force $M$ to have the form $(M_1,M_2,O_{2\times2})$ where $M_1=M_2=(1)$, and so $M$ would only contain two ones. Next, if $q=1$, then again we apply \ref{bound on l for 4} to give ${\rm rk} (H \vee K) \le 1+{\rm min}(2+\frac{2}{2}, 2+\frac{1}{2})=1+{\rm min}(3,2.5)=3.5$. Of course, this says that ${\rm rk} (H \vee K) \le 3$ as ranks must be integral and the conclusion is established.\EndProof

We now restate Theorem \ref{true for 5} from the Introduction:
\Theorem \label{true for 5 again} Suppose $M$ is a closed,
orientable, hyperbolic $3$-manifold such that $\pi_1 (M)$ is
$k$-free with $k = 5$. Then when $\lambda = \log9$, we have $r_M
\le 2$.\EndTheorem

\Proof This is a direct result of Corollary \ref{Corollary for 5} along with Theorem \ref{rM is less than k-3 again}.
\EndProof

For the final Corollary recall Definitions \ref{mathfrak X}.
\Corollary \label{cor for alts for hyp sp} Suppose $M={\bf H}^3 /
{\Gamma}$ is a closed, orientable hyperbolic $3$-manifold with $\Gamma \le {\rm Isom}_+({\bf H}^3)$ discrete, purely loxodromic and $5$-free. Then when $\lambda =\log9$, one of the following three alternatives holds: \Alternatives

\item \label{cor hyp sp alt0} $M$ contains an embedded ball of radius $(\log9)/2$,

\item \label{cor hyp sp alt1} There exists a point $P^{*}
\in {\bf H}^3$ with $P^{*} \in {\mathfrak X}_M$ such that ${\mathfrak
s}_m(P^{*})$, is equal to $\log9$, or

\item \label{cor hyp sp alt2} ${\bf H}^3 =
{\bigcup}_{C_1,C_2 \in C_{\log{9}}(\Gamma)} Z_{\log{9}}(C_1) \cap
Z_{\log{9}}(C_2)$ (i.e. ${\rm rk }G_{\widetilde P} \ge 2$ for all ${\widetilde P} \in {\bf
H}^3$), \emph{and} there exists a point ${\widetilde P} \in {\bf H}^3$ such that ${\rm
rk} \; G_{\widetilde P} = 2$.

Let $q: {\bf H}^3 \rightarrow {\bf H}^3 / \Gamma$ be the projection map. As $M= {\bf H}^3 / \Gamma$, we have $\Gamma \cong \pi_1(M)$. We may then equivalently restate (\ref{cor hyp sp alt2}) to say there exists a point $P=q({\widetilde P})$ in $M$ such that the class of all homotopically non-trivial loops of $\pi_1(M,P)$ of length $\le \log9$ is contained in a rank-$2$ subgroup of $\Gamma$. \EndAlternatives \EndCorollary

\Proof The result of Theorem \ref{true for 5 again} is that $r_M \le 2$; so the only possible values for $r_M$ are $0$, $1$ and $2$.

Case (\ref{cor hyp sp alt0}) follows
when $r_M=0$ and is the result of Proposition \ref{rank zero}, and Case (\ref{cor hyp sp alt1}) follows when $r_M=1$ and is the result of Proposition \ref{rank one}. Case (\ref{cor hyp sp alt2}) occurs when $r_M=2$ and is merely restating that definition.
\EndProof

\bibliographystyle{plain}

\end{document}